\
\magnification 1100

\font\tenmsb=msbm10
\font\sevenmsb=msbm10 at 7pt
\font\fivemsb=msbm10 at 5pt
\newfam\msbfam
\textfont\msbfam=\tenmsb
\scriptfont\msbfam=\sevenmsb
\scriptscriptfont\msbfam=\fivemsb
\def\Bbb#1{{\fam\msbfam\relax#1}}

 \def\la{\longrightarrow}
 \def\ni{\noindent}
 \def\cl{\centerline}
 
\def\d{\delta }

\def\g{\gamma}

 \def\O{{\cal O}}
 \def\Y{{\cal Y}}
 \def\X{{\cal X}}

\def\P{{\Bbb P}}
\def\Q{{\Bbb Q}}

\def\sqr#1#2{{\vcenter{\vbox{\hrule height.#2pt \hbox{\vrule width.#2pt height
#1pt \kern #1pt \vrule width.#2pt}\hrule height.#2pt}}}}

\def\prod{d_1 d_2}

\def\propb{2.5}

\overfullrule=0pt

\def\Prod{(D_1\cdot D_2)}
\def\ProdE{(D_1\cdot E)(D_2\cdot E)}
\def\ProdEE{(D_1\cdot E)(D_2\cdot E)(D_3\cdot E)}

 \def\dual{\omega _{{\cal Y}/B}}

 \def\Pic{{\rm Pic}}
\def\NS{{\rm NS}}
\def\NeSe{N\'eron-Severi \ }

\def\pull{\pi ^*}

\def\coeffL{\{ \pi ^*L\} }
\def\coeffC{\{ \pi ^*C\} }
\def\JL{J^L}
\def\GL{G^L}
\def\KL{K^L}
\def\KKL{K'^L}
\def\HL{H^L}
\def\JM{J^M}
\def\KM{K^M}
\def\KKM{K'^M}
\def\HM{H^M}
\def\Jb{J_{2,b}}
\def\Gb{G_{E,b}}
\def\Hb{H_{2,b}}

\def\Htb{H_{3,b}}
\def\HEb{H_{E,b}}
\def\Kb{K_{2,b}}
\def\Kub{K_{1,b}}
\def\KEb{K_{E,b}}
\def\Kob{K_{0,b}}
\def\KKb{K'_{2,b}}

\def\KKEb{K'_{E,b}}
\def\KKob{K'_{0,b}}

\def\Pr{{\Bbb P} ^r}
\def\Pl{{\Bbb P} ^2}

\def\Ou{{\cal O}_{\Pl }(1)}
\def\F3{{\Bbb F}_3}
\def\one{{\Bbb F}_1}
\def\two{{\Bbb F}_2}
\def\Fn{{\Bbb F}_n}
\def\E{\tilde E}
\def\ND{N(D)}
\def\NDD{N(D')}

\def\Nd{N(d)}

\def\tD{N_2(D)}

\def\tDDD{N_2(D'')}
\def\fDE{N_3(D - E)}
\def\fD{N_3(D)}
\def\rD{r_0(D)}
\def\rDu{r_0(D_1)}
\def\rDt{r_0(D_2)}

\def\VD{V(D)}
\def\TD{V_2(D)}
\def\bound{\Delta }

\def\GE{{\bf G}_E}
\def\G{{\bf G}}
\def\J{{\bf J}}
\def\K{{\bf K}}
\def\KK{{\bf K}'}
\def\H{{\bf H}}
\def\j{j(D_1,D_2)}
\def\k{k(D_1,D_2)}
\def\kk{k'(D_1,D_2)}
\def\h{h(D_1,D_2,D_3)}
\def\htwo{h(D_1,D_2)}
\def\nJ{n_J}
\def\nH{n_H}
\def\nK{n_K}
\def\nKK{n_{K'}}
\def\su{\sum _{d_1+d_2 =d}}
\def\summ{\sum _{D_1+D_2 =D}}
\def\summm{\sum _{D_1+D_2 =D-E}}
\def\summmm{\sum _{D_1+D_2 + D_3=D-E}}

\

\

\centerline{\bf{ENUMERATING RATIONAL CURVES:}}

\centerline{\bf{THE RATIONAL FIBRATION METHOD}}

\

\

\

\

\noindent {\bf Lucia Caporaso}

\noindent Mathematics department, Harvard University,

\noindent 1 Oxford st., Cambridge MA 02138, USA

\noindent caporaso@abel.harvard.edu

\

\noindent {\bf Joe Harris}

\noindent Mathematics department, Harvard University,

\noindent 1 Oxford st., Cambridge MA 02138, USA

\noindent harris@abel.harvard.edu

\

\

\

\ni \cl{{\bf Contents}}

\

{\narrower\narrower\narrower\smallskip
\noindent 1. Introduction \thinspace \dotfill \thinspace 1

\noindent 2. The rational fibration method in general  \thinspace \dotfill
\thinspace 4

\ni 3. Plane curves \thinspace \dotfill \thinspace 6

\ni 4. The general recursion for $\two$  \thinspace \dotfill \thinspace 12

\ni 5. The general recursion for $\F3$  \thinspace \dotfill \thinspace 15

\ni References  \thinspace \dotfill \thinspace 24
\smallskip}

\

\

\

\ni \cl{{\bf 1. Introduction}}

\

\ni {\it  $1.1.$ Statement of the problem.}

 Let $S$ be a smooth projective rational
surface and let
$D$ be an effective divisor on
$S$. Let $V(D) \subset |D|$ be the closure of the locus of irreducible
rational curves. For general results about the
geometry of $V(D)$,  we refer to [H] and to [CH].

If $D$ has nonnegative self-intersection and $\VD$ is nonempty the
dimension  of $\VD$ is known (cf. [K]):
$$
\rD : = \dim \VD \; = \;  -(K_S\cdot D) -1.
$$

\ni The  problem that we will study here is  to compute the degrees
$$
\ND \; := \;  \deg \VD
$$
of these varieties as subvarieties of
$|D| \cong \Pr$. Alternatively, $\ND$ is the number of irreducible rational
curves in
$|D|$ that pass through  $\rD$  general points of $S$. If $S$ is the
projective plane
$\Pl$ and $d= \deg D$, then one also uses the notation
$\Nd$ to denote the number of irreducible, rational curves of degree $d$ passing
through $3d-1$ general points.

\

\ni {\it $1.2.$ Terminology and notation.}

We will work over the
complex numbers. Throughout, the words ``surface" and ``curve" will refer to
projective varieties.

If $D$ and $D'$ are effective
divisors (or divisor classes) on a surface, we will say that $D>D'$ if
$D-D'$ is effective and
nonzero.

We will denote by
$\Fn$ the {\it Hirzebruch surface} $\Fn = \P(\O_{\P^1} \oplus
\O_{\P^1}(n))$. On each
$\Fn$ with $n \ge 1$ there exists a unique curve of negative self
intersection, which we will
denote by
$E$ and refer to as the {\it exceptional curve} on $\Fn$. We will denote by
$F$ a
fiber of the projection $\Fn \to \P^1$; the classes of $E$ and $F$ generate
the Picard group of
$\Fn$, with intersection pairing given by
$$
E^2 = -n ; \qquad (E \cdot F) = 1 \qquad {\rm and} \qquad F^2 = 0 .
$$
Another useful divisor class  is the
 class of a complementary section, that is, a section $C$ of
the $\P^1$-bundle $\Fn \to \P^1$ disjoint from $E$. Since $(C \cdot E) = 0$
and $(C \cdot F) = 1$,
we see that  $C \equiv E + nF$; so the classes $C$ and $F$ also generate the
Picard group, with intersection  numbers
$$
C^2 = n ; \qquad (C \cdot F) = 1 \qquad {\rm and} \qquad F^2 = 0 .
$$
For any positive integer $m$, we will denote by $V_m(D) \subset V(D)
\subset |D|$ the closure of the locus of
irreducible rational curves $X$ having contact of order at least $m$ with
$E$ at a smooth point of $X$. These
varieties will also be referred to as Severi varieties. We set $N_m(D) :=
\deg V_m(D)$.

\

\ni {\it $1.3.$ Methods and results.}

Until very recently, the basic enumerative problem of determining the
degrees  of
Severi varieties was unsolved even in the case of
 $\Pl$. In 1989 Ziv Ran [R] described a recursive procedure for calculating the
degrees of the Severi varieties parametrizing plane curves of any degree and
genus. Recently, M. Kontsevich  discovered a beautiful and simple
recursive
 formula in the case of rational curves on
$\Pl$ (see [KM] and [RT] for proofs).  Kontsevich's method was
based on his description of a compactified moduli space for maps of $\P^1$
into the
surface $S = \P^2$;  others (e.g., [DI], [KP] and [CM]) were able to use
the same method
to derive similar formulas in the case of other surfaces $S$ for which a
Kontsevich-style moduli space existed, such as $S = \P^1
\times
\P^1$, the ruled surface $S =
\one$ and del Pezzo surfaces.

It was our feeling that the reliance of Kontsevich's method on the
existence of a
well-behaved moduli space was not essential. We were especially interested in
whether a similar formula might be derived for the Hirzebruch surfaces $\Fn$.

In [CH], we succeeded in recasting the
Kontsevich method so as to remove the apparent dependence on the existence of a
moduli space: as we set it up, it was necessary only to understand the
degenerations of the
rational curves in the one-parameter families corresponding to general
one-dimensional linear
sections of $\VD$.
The resulting ``cross-ratio method" allowed us to derive a complete
recursion for all divisor classes on the ruled surface $S =
\two$---that is, a formula expressing $\ND$ in terms of $N(D')$ for $D' <
D$---and a closed-form
formula for certain divisor classes on the ruled surfaces $\Fn$ for any
$n$. (In fact, compactifications of the moduli space of maps $\P^1 \to S$
do exist for these
surfaces, but they contain in general many components, only one of which
parametrizes
generically irreducible rational curves and the others of which may have
strictly larger
dimension. Kontsevich's method can be carried out in these cases, as was
done by Kleiman and
Piene [KP]; but at present we do not see how to use the resulting formulas to enumerate
irreducible rational curves.)

However, we were unable to go significantly beyond this point: a
similarly derived formula in [CH] for the degrees $N(D)$ of Severi
varieties $V(D)$ on $\Fn$
expresses
$N(D)$ not solely in terms of
$N(D')$ for
$D' < D$, but also in terms of the degrees $N_k(D'')$ of the Severi
varieties $V_k(D)$
parametrizing curves with a point of $k$-fold tangency with a fixed curve
$E \subset S$.  For
example, if $S=
\F3$, then
$\ND$  is expressed as a function of
$N(D')$ and of
$\tDDD$, where
$\tDDD$ is the number of irreducible rational curves in
$|D''|$ that are simply tangent to $E$ and pass through the appropriate
number (that is,
$r_0(D'') -1$) of general points of $\F3$. A complete recursion in this
case would have required a
similar analysis of linear sections of the Severi varieties $V_2(D)$, which
in turn would have
have necessitated an analysis of Severi varieties parametrizing curves with more
complicated tangency conditions.

In the end, it seems that one way or another we need to deal with the
degrees of these
``tangential" loci as well. This difficulty led us to the discovery of a
computational technique
different from and simpler than the cross-ratio method, which we will
describe in the present
paper. It involves an analysis of the same basic object as the cross-ratio
method---that is, the
one-parameter family
$\X \to
\Gamma$ of rational curves through $\rD - 1$ general points of
$S$ and their limits---but extracts more information from it. It is  based
on a description of the
N\'eron-Severi group of a minimal desingularization  of $\X$ (we will
therefore refer to it as
the ``rational fibration method"). The main advantage of
this technique for our present purposes is that we are in fact able compute
the degrees of the
tangential loci involved; at least in all cases that we studied. It also
yields other related formulas,
such as the number of irreducible rational curves having a node at a given
general point $p \in
S$ and passing through $r_0(D) -2$ other  general points.

\

\ni {\it $1.4.$ Contents of this paper.}

In the following section we will describe the rational fibration method in
a general setting. In
the succeeding sections we  will apply it in the cases $S=\P^2$, $S=\two$
and $S=\F3$. In the
first of these cases, we obtain another (simpler) proof of Kontsevich's
formula, as well as  some
related formulas derived by Pandharipande [P]. In the second, we will
recover the general
recursion formula found originally in [CH] for degrees of Severi varieties
on $\two$. Finally, in
the last section we derive a complete set of recursions for $\F3$, the
first case for which the
cross-ratio method does not give a complete answer.

We have tried to keep this paper relatively self-contained; in particular
it should be
intelligible to a reader unfamiliar with [CH]. We will, however, have to
appeal to chapter 2
of [CH] for proofs of some of the basic assertions about the local geometry
of Severi varieties
and the families of curves they parametrize.

\

\

\ni {\cl{ \bf 2. The rational fibration method in general}}

\

\ni {\it $2.1.$ Objects and morphisms.}

As we indicated in the introduction, the rational fibration method, like
the cross-ratio method,
involves studying a suitably general one-parameter family of rational
curves. To set it up,
first let $S$ be a smooth rational surface and
$D$ an effective divisor on $S$. We will assume that $D$ has nonnegative
self-intersection
and that the Severi variety $V(D) \ne \emptyset$, so that in particular we
have $\dim V(D) =
\rD = -(K_S \cdot D) -1$. Now, choose $\rD - 1$ general points
$p_1,\ldots,p_{\rD-1} \in S$
and let $\Gamma \subset \VD$  be the closure of the locus of points $[X]
\in V(D)$
corresponding to irreducible rational  curves
$X$  passing through these points. Equivalently, if for any point $p \in S$
we let $H_p \subset
\P^r$ be the hyperplane of points  corresponding to  curves
  passing through $p$,
$\Gamma
$ will be the one-dimensional linear section of
$\VD$
$$
\Gamma \; = \; V(D) \cap \bigl(\bigcap_{i=1}^{r_0(D)-1} H_{p_i} \bigr) .
$$
Now, let $\X \subset \Gamma \times S \to \Gamma$ be the  family of curves
corresponding
to $\Gamma \subset |D|$.  Consider the  normalization $\Gamma^\nu \to \Gamma$ of
$\Gamma$, and then take the normalization
$\X^\nu$ of
$X
\times_\Gamma \Gamma^\nu$  to arrive at a family
$$
\X^\nu \; \la \; \Gamma^\nu
$$
over a smooth curve $\Gamma^\nu$, whose general fiber is isomorphic to $\P^1$.

Next, we  apply  semi-stable reduction (which we should rather call
``nodal reduction",
since our curves have genus zero): after making a base change $B \to
\Gamma^\nu$ and
blowing up the total space of the pullback family $\X^\nu
\times_{\Gamma^\nu} B$, we
arrive at a family
$f : \Y
\to B$ whose total space is smooth, whose general fiber is a smooth
rational curve and whose
special fibers are all nodal curves. In fact, a base change will turn out
to be unnecessary in
each of the three cases considered below---the minimal desingularization
$\Y$ of the total
space
$\X^\nu$  already has this property---but this is not relevant, since
even a superfluous base change will not affect the subsequent calculations.
We will denote
by $\pi : \Y \to S$ the composite map
$$
\pi \; : \; \Y \to \X^\nu \times_{\Gamma^\nu} B \to \X^\nu \to \X \subset
\Gamma \times S
\to S .
$$
Notice that $\pi$ is a generically finite map, whose degree
is  the number of irreducible
rational curves in the linear series
$|D|$ passing through the points $p_1,\ldots,p_{\rD-1}$ and $p$: that is,
the degree $N(D)$ of
the Severi variety $V(D)$.

Here is a diagram of the basic objects and morphisms we have introduced:

\vfill\eject

\

\vskip2.9in

\hskip.8in \special {picture basicdiagram}

\

\ni {\it $2.2.$ Outline of the method.}

As we indicated, our method involves calculating in the N\'eron-Severi group
of the total space $\Y$ of our family. This is motivated by a simple
observation: given any two
line bundles
$L$ and
$M$ on
$S$, we have
$$
(\pull L \cdot \pull M) \; = \; \ND (L\cdot M).
$$
Thus, in order to
derive  a formula for $\ND$,  we want to compute intersection numbers in
the \NeSe  group
of $\Y$. For example, if  $S$ is the projective plane $\Pl$ we
can take
$L=M=\Ou$. Then $(\pull L)^2=\ND$; so if we can compute
$(\pull L)^2$  we get a formula for $\ND$.

What makes it possible to perform such calculations if the fact that $f$
expresses $\Y$
 as the total space of a one-parameter family of generically smooth
rational curves, so
that to determine the class of a given divisor it is enough to know its
degree on each
component of each fiber of $f$. More precisely, the Picard group of $\Y$ will be
freely generated by the class of a fiber, the class of any section $A$ of
$f$, and the classes of
all the irreducible curves contained in fibers of $f$ and disjoint from
$A$. Moreover, in terms
of these generators the intersection pairing on $\Pic(\Y)$ is (except
possibly for the
self-intersection of $A$) easy to describe. This means two things: first,
we can express a given
divisor class as a linear combination  of these generators once we know its
degree on each
component of each reducible fiber and on
$A$; and second, having expressed two divisor classes as linear
combinations of these
generators, we can readily compute their intersection number.

The method we will apply in each case thus consists of five steps:

\

\item{$\bullet$} \thinspace First, we need to describe the reducible
fibers of $\Y \to B$; that is (given that $B$ will be in practice just the
normalization
$\Gamma^\nu$ of the base of our original family $\X \to \Gamma$), the set
of reducible
curves in the linear series $|D|$ through the points $p_1,\ldots,p_{\rD-1}$
that are limits of
irreducible rational curves through these points, and the branches of
$\Gamma$ at each one.
The characterization of such curves is straightforward in the case of
$S = \P^2$ by simple dimension-counting. In the case of $S = \Fn$ with $n
\ge 2$ it is less
obvious, since in contrast with the case of $\P^2$ most reducible curves
through the points
$p_1,\ldots,p_{\rD-1}$ whose components are all rational are not  limits of irreducible
rational such curves; the answer is worked out in [CH]. In either case the
number of such
fibers will be known inductively.   \smallskip

\item{$\bullet$} \thinspace Secondly, we need to describe the local
structure of the
family
$\X^\nu
\to
\Gamma^\nu$ near each reducible fiber; specifically, we need to know
whether $\X^\nu$ is
smooth or if we have to blow up. This likewise is straightforward in the
case of the plane,
where in fact
$\X^\nu$ is smooth. It is more interesting in the case of the Hirzebruch
surfaces $\Fn$, where
for $n \ge 3$ we see that $\X^\nu$ will indeed have singularities; again,
this is worked out in
[CH] and we will refer there for the relevant results.   \smallskip

\item{$\bullet$} \thinspace  Third, we choose a basis for the
N\'eron-Severi group of $\Y$, and calculate the intersection pairing on
these classes.
\smallskip

\item{$\bullet$} \thinspace  Fourth,
since we know the images in $S$
of the components of reducible fibers of $f : \Y \to B$, we can calculate
the degrees on all such
components of the pullback $\pi^*L$ of any line bundle $L$ on $S$; and
\smallskip

\item{$\bullet$} \thinspace  Fifth, we are able
therefore to express the intersection numbers $(\pi^*L \cdot \pi^*M)$ for
pairs $L,M \in \Pic
S$ of line bundles on $S$.

\

\ni
Evidently, the particulars of this process will depend on
$S$ and $D$; for the moment we shall just fix some notation and make some
preliminary
observations. First, for $b\in B$ we use the common notation $Y_b:=f^{-1}
(b)$ to denote the
fiber of $f$ over
$b$. The class in $\NS (\Y)$ of such a fiber is denoted by $Y$.

Secondly, recall that our family parametrizes curves through certain base
points. We pick two
of them,
$q$ and
$q'$, and we denote by $A$ and $A'$ the corresponding sections of $f$. The
following
relations are clear:

$$Y^2 = A\cdot A' = 0 \  \   {\rm and} \  \  A\cdot Y =A'\cdot Y = 1.$$
Notice also
that by symmetry $A^2 = {1 \over 2}(A-A')^2$ which will be useful to
compute the left hand
side. In fact
$A-A'$ is supported on exactly those fibers  where $q$ and $q'$ lie on different
components, the number of which we will be able to  count.

One further note: the description above of the N\'eron-Severi group of $\Y$
as generated by
the classes of $A$, $Y$ and components of reducible fibers assumes that the
base $B$ of the
family is connected, which we will not always know in practice.  This
assumption is not
essential, however: in case  $B$ has irreducible components
$B_1,\ldots,B_k$ we simply have
to replace every multiple of $Y$ in the formulas below by a suitable linear
combination of
fibers $Y_i$ lying over points of $B_i$. As the reader may verify, this
does not alter the
outcome of the subsequent calculations.

\

\

\cl{\bf 3.  Plane curves}

\

\ni Here we study the case $S=\Pl$. If $D$ and $D_i$ are divisors on $\Pl$,
we denote
their degrees respectively  by $d$ and
$d_i$. Since a divisor class in the plane is  determined by its degree, we
will introduce the  notation $\Nd :=\ND $.

We have that $\rD = 3d -1$, so we choose general points
$p_1,\ldots,p_{3d-2} \in \P^2$, let
$\Gamma \subset V(D)$ be the locus of curves in $V(D)$ containing the
points $p_i$, and
proceed as described in the preceding section. To describe the resulting
family of curves, let
$\bound$ be the locus of
$\VD
$ parametrizing
 degenerate curves (that is, curves that are reducible or have
singularities other
than nodes). Since our curve $\Gamma \subset V(D)$ will intersect $\Delta$
only at general
points of components of $\Delta$, we may apply the results of [DH] and [H]
to conclude the
following

\

{\narrower \smallskip \ni A. \thinspace  Any fiber $X_\g$ of $\X \to
\Gamma$ is either
\item{1.} an irreducible curve with exactly $\d={(d-1)(d-2) \over 2}$ nodes;
\item{2.} an irreducible curve with exactly $\d-1$ nodes and a cusp;
\item{3.} an irreducible curve with exactly $\d-2$ nodes and a tacnode;
\item{4.} an irreducible curve with exactly $\d-3$ nodes and an ordinary
triple point; or
\item{5.} a curve having exactly two irreducible components $X_1$, $X_2$, of degrees
$d_1$ and
$d_2$, with exactly ${(d_1-1)(d_1-2) \over 2}$ and ${(d_2-1)(d_2-2) \over
2}$ nodes
respectively, and intersecting transversally in $d_1d_2$ points.

\

\ni B. \thinspace In cases $1$, $3$ and $4$, the curve $\Gamma$ is smooth
at $\g$ and the
family
$\X^\nu
\to
\Gamma^\nu$ is smooth at the unique point of $\Gamma^\nu$ lying over $\g$.
In case
$2$, $\Gamma$ has a cusp at $\gamma$ but the family $\X^\nu \to
\Gamma^\nu$ is still smooth at the unique point of $\Gamma^\nu$ lying over
$\g$.

\

\ni C. \thinspace In
case $5$ the curve $\Gamma$ has $d_1d_2$ smooth branches at $\g$
(corresponding to
deformations of $X_\g$ smoothing any one of the $d_1d_2$ nodes of $X_\g$
coming from a
point
$p
\in X_1 \cap X_2$ of intersection of $X_1$ and $X_2$). At each point of
$\Gamma^\nu$
lying over $\g$ the fiber of the family
$\X^\nu
\to
\Gamma^\nu$ has two smooth rational components meeting transversally at one
point (more
precisely, it is the normalization of
$X_\g$ at the remaining
${(d_1-1)(d_1-2)
\over 2}+{(d_2-1)(d_2-2) \over 2}+d_1d_2-1 = \d$ nodes of $X_\g$), and
smooth total space.

\

\ni D. \thinspace Finally, if $X \subset \P^2$ is any curve of type 1-5
passing through the
points
$p_1,\ldots,p_{3d-2}$, then conversely $[X] \in \Gamma$; that is, $X$ is a
limit of irreducible
rational curves $X_\g$ through $p_1,\ldots,p_{3d-2}$.
\smallskip}

\

\ni We see in particular that the total space $\X^\nu$ is smooth and that
the fibers of
$\X^\nu \to \Gamma^\nu$ are all nodal, so that no further base changes or
blow-ups are
necessary; that is, we may take
$B =
\Gamma^\nu$ (as we stated earlier) and $\Y =
\X^\nu$. Note also that every reducible fiber of
$\Y$ has precisely two irreducible components, meeting transversally at one
point.

We shall call a reducible fiber of
$f$ a {\it fiber of type}
 $\J$ and we shall denote by $B_J$ the subset of points of $B$ such that the
corresponding fiber is of type $\J$, that is, reducible. (This  new piece of
terminology probably seems pointless, but it will be useful  in the
sequel). For any
$b\in B_J$, then we denote the two irreducible components of the fiber over
$b$  by
$J_{1,b}$ and $\Jb$. We shall always denote by
$J_{1,b}$ the component containing the point $q=p_1$. The picture of $\Y$ is
thus:

\vfill\eject

\

\vskip2in

\hskip1in \special {picture P2picture}

\

\ni
Let
$D_i$ be the class of
$\pi (J_{i,b})$.   We denote by $\j$ the number of all such fibers, for any
given
decomposition $D=D_1+D_2$. To determine $\j$, note first that if $X=X_1
\cup X_2 \in
\Gamma$ is any reducible curve, $D_i$ the class of $X_i$, then $X_i$ can
contain at most
$r_0(D_i)$ of the
$r_0(D)-1 = 3d-2$ points
$p_1,\ldots,p_{3d-2}$. Since
$$
r_0(D_1) + r_0(D_2) \; = \; r_0(D) -1 ,
$$
it follows that each component $X_i$ must contain exactly $r_0(D_i)$ of the
points
$p_1,\ldots,p_{3d-2}$. Thus, to specify such a curve, we have first to
choose a decomposition
of the set
$\Phi = \{p_1,\ldots,p_{3d-2}\}$ into disjoint subsets $\Phi_1$, $\Phi_2$
of cardinality
$r_0(D_1)$ and
$r_0(D_2)$ respectively, with the point $q=p_1 \in \Phi_1$; and then to
choose, for each $i$,
one of the
$N(D_i)$ curves
$X_i
\in V(D_i)$ containing $\Phi_i$. The number of such curves $X$ is thus
$$
{r_0(D)-2 \choose r_0(D_1)-1} \cdot N(D_1) \cdot N(D_2)
$$
and since we have seen there are $(D_1\cdot D_2) = d_1d_2$ points of $B$
lying over each
point $[X] \in
\Gamma$ corresponding to a curve of this type, we have
$$
\j \; = \; N(D_1)N(D_2)(D_1\cdot D_2) {\rD -2 \choose \rDu - 1}\; = \;
N(d_1)N(d_2)\prod {3d-3
\choose 3d_1 -2 }.
$$
Note that by a simple dimension count, any of the curves $X_i \in V(D_i)$
passing through
$r_0(D_i)$ of the
points
$p_1,\ldots,p_{3d-2}$ will be irreducible and nodal. By a standard further
argument as in
Lemma 2.1 of [CH], we see that any pair $X_1$, $X_2$ of such curves will
intersect
transversally, so that the union $X = X_1 \cup X_2$ will indeed be a curve
as described in
(5) above.

This completes the first two steps in the general method. Next, we give a basis for the
N\'eron-Severi group
$\NS (\Y )$. We now choose as a system of generators for $\NS (\Y )$ the
class $A$ of the
section of $f :
\Y \to B$ coming from the base point $q = p_1$ of our family; the class $Y$
of a fiber of the
map $f : \Y \to B$, and the classes $\{J_{2,b}\}_{b \in B_J}$. Most of the
pairwise intersection
numbers of these classes are readily given: we clearly have
$$
\eqalign{(A\cdot Y) \;  &= \; 1 \cr
Y^2  \;  &= \; 0 \cr
(A \cdot J_{2,b})   \;  &= \; 0  \quad \forall \; b \cr
(Y \cdot J_{2,b})   \;  &= \; 0  \quad \forall \; b  \cr
(J_{2,b} \cdot J_{2,b'})   \;  &= \; 0 \quad \forall \; b \ne b'\, ; \quad
{\rm and} \cr
(J_{2,b} \cdot J_{2,b})   \;  &= \; -1 \quad \forall \; b \, . \cr}
$$
In fact, there is only one intersection number that is not evident:  $A^2$.
To compute it we  choose a base point $q' \ne q$, so that  $q'$ determines
a second section $A'$ of
$\Y \la B$ disjoint from $A$. Since the base points $p_1,\ldots,p_{3d-2}$
of our family are
general points in the plane, by symmetry we have $A'^2 = A^2$; hence we can
write
$$
2A^2 \; = \; (A-A')^2.
$$
To compute the right hand side, let
$$
S_J=\{ b\in B_J \  {\rm{such \  that }}\  q'\in \pi(\Jb) \}
$$
be the collection of points in $B$ over which the sections $A$ and $A'$
meet different
components of the fiber; let $n_J = |S_J|$ be the cardinality of $S_J$.
For every $b \notin  S_J$, $A$ and $A'$ have the same intersection number
with each
component of the fiber $Y_b$. For $b \in S_J$, on the other hand, we have
$(A \cdot J_{1,b}) =
1$ and
$(A
\cdot J_{2,b}) = 0$, while $(A' \cdot J_{1,b}) = 0$ and $(A' \cdot
J_{2,b}) = 1$. It follows that the classes
$$
A \qquad {\rm and} \qquad A' \; - \; \sum_{b \in S_J} J_{2,b}
$$
have the same intersection number with every component of every fiber of
$\Y \la B$, and so
must differ by a multiple of the class $Y$ of a fiber: that is,
$$
A-A' \; = \;  - \sum _{b\in S_J} \Jb + nY
$$
for some integer $n$. In fact, $n$ must be equal to $n_J/2$ by symmetry,
but that is
irrelevant in any case: squaring both sides, we find that
$$
(A-A')^2 \; = \;   \sum _{b\in S_J} \Jb^2 \; = \; -n_J
$$
and hence
$$
A^2 = -{\nJ \over 2} .
$$
Thus, it remains only to determine the number $n_J$ of reducible fibers of
$\Y \la
B$ lying over curves in our original family in which the points $q$ and
$q'$ lie in different
components. We can do this in exactly the same way as we determined the
total number
of reducible fibers: the only difference is that now we want to count only
decompositions
$\Phi = \Phi_1 \cup \Phi_2$ in which $q=p_1 \in \Phi_1$ and $q' = p_2 \in
\Phi_2$. We thus
replace the binomial coefficient ${\rD -2 \choose \rDu - 1}$ in the formula
for $j(D_1,D_2)$
above with ${\rD -3 \choose \rDu - 1}$ and sum over all pairs $D_1$, $D_2$
with $D_1 + D_2
= D$ to obtain
$$
\nJ = \summ N(D_1)N(D_2)(D_1\cdot D_2) {\rD -3 \choose \rDu -1 } .
$$
This completes the third step of the process.

Now let
$L\in
\Pic (\Pl )$ be any line bundle on the plane, and write the class of its
pullback to $\Y$ as a
general linear combination of our chosen generators
$$
\pull L \; \equiv \;  c_aA \, + \, c_yY \; + \sum_{b \in B_J} c_bJ_{2,b} .
$$
We will denote by $\JL$ the third term on the right, that is, we set
$$
\JL := \sum _{b\in B_J} c_b \Jb \, ;
$$
this notation is not immediately useful, but will become so in the
succeeding calculations.

We now intersect both sides of the above equivalence with each of our
chosen generators of
$\NS(\Y)$ to determine the coefficients $c_a$, $c_y$ and $c_b$. First, by
intersecting  both
sides with
$Y$ we find that
$$c_a \; = \; (\pi^*L \cdot Y)  \; = \; (L \cdot \pi_*Y)  \; = \; (L\cdot
D) \, .
$$
Next we intersect with
$A$: we have
$$
(\pi^*L \cdot A)  \; = \; (L \cdot \pi_*A)  \; = \; 0
$$
since $\pi$ is constant on the curve $A$; and hence
$$
c_y \; = \; -A^2(L\cdot D)  \; = \;  {n_J \over 2}(L \cdot D)\, .
$$
Finally, to determine $c_b$ we naturally intersect both sides with the
class of $J_{2,b}$; we
find that
$$
c_b \; = \; -(\pi^*L \cdot J_{2,b})  \; = \; -(L \cdot \pi_*J_{2,b})  \; = \; -(L \cdot D_2) \, .
$$
Thus, in sum,
$$
\pi^*L \; \equiv \; (L\cdot D)A  + {n_J \over 2}(L \cdot D)Y - \sum_{b \in
B_J} (L \cdot
D_2)J_{2,b} \, .
$$

For the final step in the process, we evaluate the self-intersection of
$\pi^*L$: we find
$$
\eqalign{(\pi^*L)^2 \; &= \; {n_J \over 2}(L \cdot D)^2 - \sum_{b \in B_J}
(L \cdot
D_2)^2 \cr
&= \; \sum_{D_1+D_2=D} \Bigl[ {1 \over 2}N(D_1)N(D_2)(D_1\cdot D_2) {\rD -3
\choose \rDu -1
}(L
\cdot D)^2 \cr
&\qquad \qquad \qquad\qquad - N(D_1)N(D_2)(D_1\cdot D_2) {\rD -2 \choose
\rDu - 1}(L
\cdot D_2)^2 \Bigr] \cr} \, .
$$
Applying this in case $L = \Ou$, and recalling that $(\pi^*\Ou)^2 = N(d)$,
we have
$$
N(d) \; = \; \su N(d_1)N(d_2)\left[{d^2\prod
\over 2}{ 3d-4 \choose 3d_1 -2} -   d_1^2d_2^2{ 3d-3 \choose 3d_1 -2
}\right] \, ;
$$
and expanding out $d^2 = (d_1+d_2)^2$ and using the symmetry with respect
to $d_1$ and
$d_2$  we get the well known recursive formula of Kontsevich
$$
\Nd = \su N(d_1)N(d_2)\left[   d_1^2d_2^2{3d-4 \choose 3d_1 -2  } -\prod
^3{3d-4 \choose
3d_1 -3}\right] \, .
$$

\

\ni \underbar{Remark}. In parts (A) and (B) of the statement  of results
quoted from
[H] and [DH], we describe completely all curves $X_\g$ in the family $\X
\to \Gamma$ having
other than
$\d$ nodes, and the local geometry of $\X^\nu \to \Gamma^\nu$ along each.
This is in fact
necessary to describe the N\'eron-Severi group of $\Y$, since even in those
cases where
a fiber is irreducible it is a priori possible that $\X^\nu$ will be
singular at along such a fiber,
giving rise to a reducible fiber of $\Y \to B$. Looking back over the
preceding calculation,
though, we see that even if this did happen,  it would not affect the
outcome of the
calculation, as long as the non-nodal singularities did not occur at  base
points of the family:
while the resolution of the singularities of $\X^\nu$ would create
additional curves on
$\Y$ independent in $\NS(\Y)$, the sections $A$ and $A'$ and any line
bundle pulled back
via
$\pi$ from
$\P^2$, would all have degree 0 on these curves, and so the relations of
linear equivalence
above would still hold.

Thus it was only necessary to observe that every curve  singular at a base
point is an
irreducible curve with $\d$ nodes. Since this statement will also hold for
the families of curves
on $\Fn$ that we will be considering in the following two sections, we will
in the sequel omit
the  description of the fibers $\X \to \Gamma$ other than reducible ones.

\

As another application, we give a formula for the number $N_2(d)$ of
plane, irreducible, rational curves
$X
\subset \P^2$ of degree $d$ passing through $3d-2$ given general points and
tangent to a given general line $\ell$ in the plane. Equivalently, this is
the degree of the
subvariety
$\TD$ of
$\VD$ defined as the closure of the locus of irreducible rational curves
that are tangent to
$\ell$ in
$\Pl$ at a smooth point of $[X]$ (notice that $\TD$ has codimension $1$ in
$\VD$). To
calculate this number, let
${\tilde L} = \pi^{-1}(\ell) \subset \Y$ be the preimage of $\ell$ under
$\pi$. Then
$\tilde L$ is an irreducible smooth curve, and the morphism $f : \Y \to B$
restricts to a
finite morphism $\tilde f : \tilde L \to B$ of degree $d$ on $\tilde L$.
Moreover, the set of
fibers
$Y_b
\subset
\Y$ of
$f$ tangent to $\tilde L$---that is, such that the intersection $Y_b \cap
\tilde L$ has
cardinality strictly less than $d$---corresponds to the set of curves $X_b$
in our original
family $\X \to \Gamma$ tangent to $\ell$. Thus,
$\tD$  is equal to the degree of the ramification divisor of the morphism
$\tilde f$.

Now, using the
adjunction formula, this degree is given by
$$
N_2(d) \; = \;  (\pull \Ou )^2 + (\dual \cdot \pull \Ou ).
$$
where $\dual $ is the relative dualizing sheaf of the family. Since we have
already calculated
the class of $\pi^*\Ou$ above, it remains only to determine the class of
$\dual$ in similar
terms, and then we will be able to evaluate this expression. We do this as
for $\pi^*\Ou$: we
first express $\dual$ as a linear combination of the generators:
$$
\dual \; \equiv \; e_aA \, + \, e_yY \; + \sum_{b \in B_J} e_bJ_{2,b}
$$
and then intersect both sides with the generators of $\NS(\Y)$ to determine
the coefficients.
First, intersecting with $Y$, we find that
$$
e_a \; = \; (\dual \cdot Y) \; = \; -2
$$
and then intersecting with $A$ and using the fact that $(\dual \cdot A) = - A^2$
we find that
$$
e_y \; = \; A^2 \, .
$$
Finally, we have $(\dual \cdot J_{2,b}) = -1, \; \forall b \in B_J$, and it
follows that the
coefficients
$e_b$ are all 1. Thus, in sum,
$$
\dual \; = -2A+A^2Y+\sum _{b\in B_J}\Jb \, .
$$
We finally obtain a formula first found by Pandharipande [P]:
$$ \eqalign{N_2(d)  \; &= \;  N(d) + \su N(d_1)N(d_2)   d_1d_2^2{3d-4
\choose 3d_1 -  3} \cr
&=\su N(d_1)N(d_2)\prod \left[   d_1d_2{3d-4 \choose 3d_1 -2  } -(d_2 ^2
-d_2){3d-4
\choose 3d_1 -3}\right] \cr} \, .
$$

This technique can also be used to recover another formula of
Pandharipande, for the degree
of the closure of the locus of irreducible rational curves of degree $d$
having a cusp. To
obtain this, we simply apply Porteous' formula to the differential
$$
d(f \times \pi) \; : \; T_\Y \to (f \times \pi)^*T_{B \times \P^2}
$$
of the map $f \times \pi : \Y \to B \times \P^2$; the classes on $\Y$
involved have already
been calculated. It should also be possible to determine in similar fashion
the degrees on
$\Gamma$ of all the divisor classes introduced in [DH], and in particular
obtain formulas for
the number of irreducible rational curves through $3d-2$ points and having
a tacnode, or the
number of irreducible rational curves through $3d-2$ points and having a
triple point, etc. At
this point, however, we conclude our study of the plane and turn to the
Hirzebruch surfaces.

\

\

\cl{\bf 3. The general recursion for $\two$}

\

\ni Let now $S=\two$. Let $C$, $E$ and $F$ be the curves on $\two$
described in Section 1.
 We apply our method exactly as before: for any effective divisor class $D$
on $S$
with $V(D) \ne \emptyset$, we choose $r_0(D)-1$ general points
$p_1,\ldots,p_{r_0(D)-1}$ of
$S$ and consider the family $\X \to \Gamma$ of curves $X \in V(D)$ passing
through the
points $p_i$; we let $\X^\nu \to \Gamma^\nu$ and $\Y \to B$ be derived from
this family as
in the general set-up.

We first describe the various types of reducible fibers that our family $\Y
\la B$ has.
The following analysis is based on Propositions 2.1, 2.5, 2.6  and 2.7  of
[CH]. In particular, the various types of degenerations can be classified as an
application of Proposition 2.5, and the singularities of $\Gamma$ at the
points $\g \in
\Gamma$ corresponding to each as an application of Proposition 2.6. Moreover,
Proposition 2.7 assures us that, just as in the case of
$\Pl$, the normalization $\X^\nu$ is smooth and so the total space $\Y$
coincides with ${\cal
X}^{\nu }$. In particular, we see that no irreducible component of any
fiber of $\Y$ is mapped
to a point by
$\pi$.

\ni With that said, we have the following classification of reducible
fibers of $\Y \to B$:

{\narrower \smallskip \ni Type $\J$.  Fibers having two smooth irreducible
components $J_1$
and
$J_2$, meeting transversally at one point, such that $\pi (J_i )= D_i$ with
$D_i >0$  and not
equal to $  E$.  We will always assume that $q\in J_1$. For any decomposition
$D=D_1+D_2$ we have that the number $\j$  of fibers of type $\J$ such that $\pi
(J_i)=D_i$ is
$$
\j =N(D_1)N(D_2)(D_1\cdot D_2) {\rD -2 \choose \rDu -1 }.
$$
The factor $\Prod$ appears because, just as in the case of $\Pl$, if $[X]\in
\Gamma $ corresponds to a curve of type $\J$, in the normalization map
$\nu :B\la \Gamma$ the fiber over $[X]$ contains exactly  $\Prod $ points
(here we are using
Proposition 2.6  of [CH]).

We let $B_J$ be the subset of points $b$ in $B$ whose fiber $X_b$ is a
curve of type
$\J$.

\

\ni Type $\G$.  Fibers having two smooth irreducible components $G_1$ and $G_E$,
meeting transversally at one point, such that $\pi (G_E) = E$ and
$\pi(G_1)$ is simply
tangent to $E$. Clearly
$q\in G_1$.  The total number of such fibers will not matter in the subsequent calculation.

We let $B_G$ be the subset of points $b$ in $B$ whose fiber $X_b$ is a
curve of type
$\G$.

\

\ni Type $\H$.  Fibers having three irreducible components $H_1$, $H_2$,
$H_E$, such
that $\pi (H_E)=E$ and
$\pi (H_i)=D_i$, with $D_i > 0$ and
$D_1 + D_2 = D-E$ (again, we will choose the labelling so that $q \in H_1$
always). By
Proposition 2.6 of [CH], if
$[X]\in
\Gamma$ is a point corresponding to this type of curve, then the fiber of
$B$ over $\Gamma
$ contains exactly
$\ProdE$ points. Hence the total number of fibers of type $\H$ that
correspond to a given
decomposition
$D= D_1 + D_2 + E$ is
$$
\htwo = \summ N(D_1 ) N(D_2) \ProdE {\rD -2 \choose \rDu -1}
$$ And just like for the other types, we define $B_H$ to be the subset of
points of
$B$ parametrizing curves of type $\H$. \smallskip}

\

\ni The picture of $\Y \la B$ thus looks like this:

\vfill\eject

\

\vskip2.5in

\hskip.5in \special {picture F2picture}

\

\

\ni Now we choose the following set of generators for the N\'eronSeveri
group of $\Y$

$$
\{ A,\  Y\} \cup  \{ \Jb\}_{b\in B_J}
\cup  \{ \Gb\}_{b\in B_G}
\cup  \{ \Hb , \  \HEb\}_{b\in B_H }.
$$

The following relations are obvious
$$
\Gb ^2 = \Jb ^2 = \Hb ^2 = -1; \  \  \HEb ^2 = -2
$$
and the intersection number of $A$ with any generator other than $A$ and $Y$ is
zero.

 We compute $A^2$ by the same
argument as used in the preceding section. We define
$$S_H=\{ b\in B_H \  {\rm{such \  that }}\  q'\in \Hb \}$$ and we let $\nH=
|S_J|$.
Then  we have
$$
\nH = \summm N(D_1)N(D_2)\ProdE { \rD - 3 \choose \rDu - 1}
$$ Similarly, we let $S_J$ and $\nJ$ to be defined exactly as in the
preceding section,
and notice that the value of $\nJ$ is expressed by the same formula  that
we had in the
plane. We obtain
$$ A-A' = - \sum _{b\in S_J} \Jb - \sum _{b\in S_H}(\HEb +2 \Hb) + nY
$$ for some integer $n$.  Hence
$$A^2 = -{\nJ +2 \nH \over 2}.$$

Let now $L \in \Pic (\two )$. We want to compute the coefficients of $\pull
L$ as a linear
combination of the chosen  generators of $\NS(\Y)$. The number of
generators being quite
large, it is now convenient to use the following notation: if
$W$ is any of the chosen generators, we shall denote by
$\coeffL _W$ the coefficient of $\pull L$ with respect to $W$. We shall
then write
$$
\pull L = \coeffL _A A + \coeffL _Y Y + \GL + \JL + \HL
$$ where $\JL$ is defined just as in the preceding section, and similarly
$$
\GL = \sum \coeffL _{\Gb }\Gb \   \  {\rm and \  } \
\HL = \sum (\coeffL _{\Hb }\Hb + \coeffL _{\HEb }\HEb ).
$$
We could now easily compute all the missing numbers in term of intersection
numbers
on $\two$; only we don't really
need it. All we need is the expression for $\pull C$; in fact we shall obtain a
formula for
$\ND$ by using the fact that $(\pull C)^2 = 2 \ND$. The following numbers
are obtained in a
straightforward way, just as in the case of
$\Pl$.

$$
\eqalign{ &\coeffC _A = (C\cdot D) = a.\cr
&\coeffC _Y = -(C\cdot D)A^2 = {a(\nJ + 2\nH) \over 2}.\cr
&\coeffC _{G_E} = 0\quad   {\rm
for \  any  \  generator  \  of \  type \  } \GE.\cr
&\coeffC _{\Jb } = -(C\cdot \pi (\Jb ) )\cr}
$$ And for any curve $H=H_1+H_2+H_E$ of type $\H$ such that $\pi (H_i
)=D_i$ we have
$$
\coeffC _{H_E} = - (C\cdot D_2) \  {\rm and \  }
\coeffC _{H_2} = - 2(C\cdot D_2)
$$

In conclusion, we get the same recursive formula that we obtained in [CH]:

\proclaim Theorem. For any effective divisor $D \ne E$ on $\two$ with $V(D) \ne
\emptyset$,
 $$
\eqalign{
&N(D) = {1\over 2}\summ N(D_1 )N(D_2)\Prod \left[{ \rD -3 \choose \rDu -1} (D_1
\cdot C)(D_2 \cdot C) - { \rD -3 \choose \rDu -2} (D_2 \cdot C)^2\right] \cr
&+\!\!\!\!\!\summm N(D_1 )N(D_2)\ProdE \left[{ \rD -3 \choose \rDu -1} (D_1
\cdot C)(D_2
\cdot C) - { \rD -3 \choose \rDu -2} (D_2 \cdot C)^2\right]\cr}
$$

\

\

\cl{\bf 4. The general recursion for $\F3$}

\

\ni Let now $S=\F3$. and let $C$, $E$ and $F$ be  as in $1.2$. Let $D$ be
an effective divisor
class
$D$ on $S$ with $V(D) \ne \emptyset$. We also introduce in this case two
additional
subvarieties of the linear series $|D|$:
 the subvariety $\TD$ of $\VD$ defined to be the closure of the locus of irreducible
rational curves $X$ tangent to $E$ at a smooth point of $X$;  and the
closure $F(D)$  of the
subvariety of
$\VD$ parametrizing irreducible curves having a smooth point of
intersection multiplicity 3
with
$E$.   Their degrees will be denoted by $\tD$ and $\fD$ respectively.

Now we proceed as before: we choose
$r_0(D)-1$ general points
$p_1,\ldots,p_{r_0(D)-1}$ of
$S$ and consider the family $\X \to \Gamma$ of curves $X \in V(D)$ passing
through the
points $p_i$; we let $\X^\nu \to \Gamma^\nu$ and $\Y \to B$ be derived from
this family as
in the general set-up. Our method will again provide us with a recursive
formula for the
degree of $\VD$, but there will be now an important difference: the
recursion will involve as
well the  degrees $\tD$ and $\fD$ of the varieties $\TD$ and $F(D)$. More
precisely, we are
going to obtain three formulas:

\

{\narrower \smallskip \ni (a) A  formula expressing $\ND$ in terms of
$\NDD$ and
$\tDDD$, where
$D'< D$ and $D''<D - E$;

\

\ni (b) A  formula expressing $\tD$ in terms of  $\ND $,
$\NDD$, $\tDDD$ and $\fDE$, where $D'< D$ and $D''<D -E$; and

\

\ni (c) A  formula expressing $\fDE$ in terms of  $\NDD$ and $\tDDD$, where
$D'< D$ and $D''<D - E$. \smallskip}

\

 We now   describe the various reducible fibers of  $\Y \la B$. Again we use
the results of [CH], in particular, Proposition \propb  and 2.7  for the
geometry
of the normalization map $B\la
\Gamma$ and of the total space ${\cal X}^{\nu }$.  By 2.7  we have that  ${\cal
X}^{\nu }$ is smooth at points lying on fibers corresponding tp types $\J$,
$\G$ and
$\H$ below; in other words, no irreducible component of a fiber  belonging
to one  of
these types is mapped to a point of $\F3$.

\

{\narrower \smallskip \ni Type $\J$.  (This is the exact analog of the type
$\J$ for $\two$)
 Fibers having two smooth irreducible components
$J_1$ and
$J_2$, meeting transversally at one point, such that $\pi (J_i )= D_i$ with
$D_i >0$
and not equal to $  E$.  We will always assume that $q\in J_1$. For any
decomposition
$D=D_1+D_2$ we have that the number $\j$  of fibers of type $\J$ such that $\pi
(J_i)=D_i$ is
$$
\j =N(D_1)N(D_2)(D_1\cdot D_2) {\rD -2 \choose \rDu -1 }.
$$ We have the coefficient $\Prod$ because, just as in the case of $\Pl$,
if $[X]\in
\Gamma $ corresponds to a curve of type $\J$, in the normalization map
$\nu :B\la \Gamma$ the fiber over $[X]$ contains exactly  $\Prod $ points.

We let $B_J$ be the subset of points $b$ in $B$ whose fiber $X_b$ is a
curve of type
$\J$.

\

\ni Type $\G$.  Fibers having two smooth irreducible components $G_1$ and $G_E$,
meeting transversally at one point, such that $\pi (G_E) = E$ and $\pi
(G_1)$ has a smooth point of contact of order 3 with $E$. Clearly
$q\in G_1$.  The total number of such fibers is $\fD$.

We let $B_G$ be the subset of points $b$ in $B$ whose fiber $X_b$ is a
curve of type
$\G$.

\

\ni Type $\K$. Fibers having four irreducible components $K_1$, $K_E$,
$K_0$ and $K_2$
forming a chain in the given order, that is $K_1\cap K_E = K_0\cap
K_E=K_0\cap K_2 =
1$  so that $K_1 ^2 = K_2 ^2 = -1$ and $\  K_E^2 = K_0^2 = -2$. As usual,
we have that
$q\in \pi (K_1)$. Moreover, $\pi (K_E) = E$ and $\pi (K_1)$ is tangent to
$E$; $\pi
(K_0)$ is a point of $E$ (namely, the point $E\cap \pi (K_2)$), in fact the
exceptional curve $K_0$  arises from the fact that the surface ${\cal
X}^{\nu }$ is
singular at the point corresponding to $E\cap \pi (K_2)$ (cf. Proposition
2.7  in
[CH]). Let $B_K$ be the subset of $B$ whose corresponding fiber is a curve
of type
$\K$. Finally, for any given decomposition $D=E+D_1 + D_2$   we have that
the number $\k$
of corresponding fibers of type $\K$ is given by
$$
\k =  N_2(D_1)N(D_2)(E\cdot D_2) {\rD -2 \choose \rDu -2 }
$$

\

\ni Type $\KK$. These are just like the fibers of type $\K$ with the only
difference
that the point $q$ belongs to the curve that is not tangent to $E$, that
is, we  have
now $\pi(K_2)$ tangent to $E$. We denote the irreducible components of such a fiber
$K'_1$,
$K'_0$,
$K'_E$
 and
$K'_2$, forming a chain in the given order, so that
$$(K'_1\cdot K'_0) =( K'_0\cdot
K'_E)=(K'_E\cdot  K_2) = 1
$$
and
$$
(K'_1 )^2 = (K'_2 )^2 = -1 \quad {\rm and} \quad  (K'_E)^2 =(
K'_0)^2 = -2.
$$
  Moreover, $\pi (K'_E) = E$  and $\pi (K'_0)$ is a point of $E$ (namely,
the point $E\cap \pi (K'_1)$). We define as usual $B_{K'}$ to be the subset
of $B$ whose
corresponding fiber is a curve of type $\KK$. Finally we see that the
number $\kk$ of such
fibers is

$$
\kk =  N(D_1)N_2(D_2)(D_1\cdot E) {\rD -2 \choose \rDu -1 }
$$

\

\ni Type $\H$.  Fibers having four irreducible components $H_1$, $H_2$,
$H_3$ and
$H_E$, such that $\pi (H_E)=E$ and
$\pi (H_i)=D_i$, with $D_i > 0$ and
$D_1 + D_2  +D_3= D-E$. In Proposition 2.6 of [CH] we proved that if
$[X]\in \Gamma$ is a point
corresponding to this type of curve, then the fiber of $B$ over $\Gamma $
contains
exactly $\ProdEE $ points. Hence the total number of fibers of type $\H$ that
correspond to a given decomposition
$D= D_1 + D_2 + D_3 +E$ is given by

$$
\h =  N(D_1)N(D_2)N(D_3)(D_1\cdot E)(D_2\cdot E)(D_3\cdot E) {\rD -2
\choose \rDu -2,  \
r_0(D_2)}
$$ And  as usual, we define $B_H$ to be the subset of points of $B$
corresponding to
curves of type $\H$. \smallskip}

\

Here is a picture displaying the various types of reducible fibers in our
family:

\vskip3in

\hskip.3in \special {picture F3picture}

\

\ni Now we choose the following set of generators for the N\'eron-Severi
group of $\Y$:

$$
\{ A,\  Y\} \cup  \{ \Jb\}_{b\in B_J}
\cup  \{ \Gb\}_{b\in B_G} \cup  \{  \Kb , \  \KEb , \Kob\}_{b\in B_K }
\cup
$$
$$
\cup  \{  \KKb , \  \KKEb , \KKob\}_{b\in B_{K'} }
\cup  \{ \Htb , \  \Hb , \  \HEb\}_{b\in B_H }.
$$
The following relations are obvious:
$$
\eqalign{\Gb ^2 = \Jb ^2 = \Kb ^2 = (\KKb ) ^2 = \Kub ^2 =\Hb ^2 = -1; \cr
\KEb ^2 = \Kob ^2
= ( \KKEb ) ^2 = (\KKob )^2 = -2; \cr \HEb ^2 = -3 \cr}
$$
and the intersection number of $A$ with any generator other than $A$ and $Y$ is
zero.

\

It will also be convenient to have a symbol denoting the class in $\NS (\Y
)$ of all
generators of the same type. Therefore we introduce the classes
$$ {\cal G}_E := \sum _{b\in B_G} \Gb, \   \  \  {\cal J}_2 := \sum _{b\in
B_J} \Jb,
$$

$$ {\cal K}_E := \sum _{b\in B_K} \KEb, \   \  \   {\cal K}_0 := \sum
_{b\in B_K}
\Kob, \  \  \   {\cal K}_2 := \sum _{b\in B_K} \Kb
$$

$$ {\cal K}'_E := \sum _{b\in B_{K'}} \KKEb, \   \  \   {\cal K}'_0 := \sum
_{b\in
B_{K'}} \KKob, \  \  \   {\cal K}'_2 := \sum _{b\in B_{K'}} \KKb
$$ and
$$ {\cal H}_E := \sum _{b\in B_H} \HEb, \   \  \   {\cal H}_3 := \sum
_{b\in B_H}
\Htb, \  \  \   {\cal H}_2 := \sum _{b\in B_H} \Hb
$$ We now use this notation immediately to write the class of the relative
dualizing
sheaf
$\dual $ of the family $f:\Y
\la B$. We have
$$
\dual  = -2A + A^2 Y + {\cal G}_E + {\cal J}_2 +  {\cal K}_E +  2{\cal K}_0
+ 3{\cal
K}_2 + {\cal K}'_0 +  2{\cal K}'_E+ 3{\cal K}'_2 +  {\cal H}_E +  2{\cal
H}_3 + 2{\cal
H}_2
$$

Now we compute $A^2$.  Let $q'$ be a
base point different from $q$, then  $q'$ determines a section $A'$ such that
 $(A\cdot A')=0$.
As before, we get
$
2A^2=(A-A')^2
$
and  we can compute the right hand side by expressing the difference $A -
A'$ as a
linear combination of components of fibers. So, let
$\nJ$ be the number of fibers of type $\J$ such that $q'$ lies on a different
component than $q$. We have

$$
\nJ = \summ N(D_1)N(D_2)(D_1\cdot D_2) {\rD -3 \choose \rDu -1 }
$$

Now let $\nK$ be the number of fibers of type $\K$ such that $q'$ lies on a
different
component than $q$; we have

$$
\nK = \summm N_2(D_1)N(D_2)(E\cdot D_2) {\rD -3 \choose \rDu -2 }.
$$
We define $\nKK $ analogously, and obtain
$$
\nKK = \summm N(D_1)N_2(D_2)(D_1\cdot E) {\rD -3 \choose r_0(D_2) -2 };
$$
and similarly $\nH$, for which we have
$$
\nH = \summmm  N(D_1)N(D_2)N(D_3)(D_1\cdot E) (D_2\cdot E) (D_3\cdot E) {\rD -3
\choose \rDu -1, \   r_0(D_2) -1}
$$
(where we will denote by ${n \choose a,\;b}$ the multinomial
$n!/a!b!(n-a-b)!$).

\

Let $S_J$ be
the subset of
$B$ consisting of those points
$b$ such that
$Y_b$ is a fiber of type
$\J$ for which $q$ and $q'$ lie on different components. Obviously $S_J$
contains $\nJ$
points. If $b\in S_J$ we write $Y_b = J_{1,b}+J_{2,b}$ (hence $q\in J_{1,b}$ and
$q' \in J_{2,b}$). In a completely analogous fashion we define $S_K$,
$S_{K'}$ and
$S_H$. If $b\in S_K$ then we write
$ Y_b=K_{1,b}+ K_{E,b}+K_{0,b}+K_{2,b}$ and similarly if $b$ is in $S_{K'}$
or $S_H$.
We therefore have that if
$b\in S_K$ (respectively, $b\in S_{K'}$  and $b\in S_H$), then $q'$ lies on
$K_{2,b}$
(respectively on $K'_{2,b}$ and $H_{2,b}$).  Now we have
$$
\eqalign{ A'-A=&\sum_{b\in S_J}J_{1,b} + \sum_{b\in
S_K}(3K_{1,b}+2K_{E,b}+K_{0,b}) +
\sum_{b\in S_{K'}}(3K'_{1,b}+2K'_{0,b}+K'_{E,b})+\cr &+ \sum_{b\in
S_H}(H_{1,b}-H_{2,b}) +nY.\cr}
$$
where $n$ is some integer that is irrelevant for our computation. Finally
we obtain
$$ A^2 = {-\nJ - 2\nH - 6\nK \over 2}
$$

\

Now, for any $L\in \Pic (\F3 )$, we  have
$$
\coeffL _A = (L\cdot D)
$$

$$
\coeffL _Y = -(L\cdot D)A^2
$$
and
$$
\coeffL _{\Gb } = -(L\cdot E)
$$
for any $b$ in $B_G$.
 These are obtained, in the given order, from the products
 $(\pull L \cdot Y) = (L\cdot D)\  $, $(\pull L \cdot A) = 0\  $  and
$(\pull L \cdot
\Gb ) = (L\cdot E)\  $

Let us fix a fiber of type $\J$ which we write as $J_{1,b} +\Jb$ as usual;
let $D_2$
be the class in $\F3$ of $\pi (\Jb )$. From the product $(\pull L \cdot\Jb)
= (L\cdot D_2)\
$ we see that
$$
\coeffL _{\Jb } = -(L\cdot D_2)
$$

Fix now a fiber of type $\K$, which we shall write as  $\Kub +\KEb +\Kob +
\Kb$, such
that the image in $\F3$  has corresponding divisor classes $D_1$ for $\Kub$
and $D_2$
for $\Kb$. The relation $(\pull L \cdot \Kub ) = (L\cdot D_1)\  $ implies

$$
\coeffL _{\KEb } = (L\cdot D_1) -(L\cdot D);
$$
the above formula together with $(\pull L \cdot
\KEb ) = (L\cdot E)\  $ gives
$$
\coeffL _{\Kob } = 2(L\cdot D_1) -2(L\cdot D)+(L\cdot E);
$$
 and the two previous formulas combined with $(\pull L \cdot K_0) = 0\  $ gives
$$
\coeffL _{\Kb } = 3(L\cdot D_1) -3(L\cdot D)+2(L\cdot E).
$$

With completely analogous notation and procedure, for a fixed fiber of type
$\KK$ we
have

$$
\coeffL _{\KKob } = (L\cdot D_1) -(L\cdot D)
$$

$$
\coeffL _{\KKEb } = 2(L\cdot D_1) -2(L\cdot D)
$$

$$
\coeffL _{\KKb } = 3(L\cdot D_1) -3(L\cdot D)+(L\cdot E)
$$ Finally,  fix a fiber of type $\H$ such that the class in $\F3$
corresponding to
$H_{i,b}$ is $D_i$; again the same procedure yields

$$
\coeffL _{\HEb } = (L\cdot D_1) -(L\cdot D)
$$

$$
\coeffL _{H_{3,b}} = (L\cdot D_1) -(L\cdot D)-(L\cdot D_3)
$$
and
$$
\coeffL _{\Hb } = (L\cdot D_1) -(L\cdot D)-(L\cdot D_2) .
$$

\ni We shall also use the following short notation:

$$
\pull L =  (L\cdot D)A + (L\cdot D)( {-\nJ - \nH - 6\nK \over 2})Y-(L\cdot
E){\cal
G}_E+
\JL + \KL + \KKL + \HL \eqno(*).
$$

\

Now we want to compute the intersection product on $\Y$ of the pull-back of two
line bundles
$L$ and $M$ on $\F3$. We easily have

$$
\eqalign{ (\pull L \cdot \pull M)=&-(L\cdot D)(M\cdot D)A^2-(L\cdot E)(M\cdot
E)f(D-E)\cr &+(\JL \cdot \JM)+(\KL \cdot \KM)+(\KKL \cdot \KKM)+(\HL \cdot
\HM).\cr}
$$
 And now a completely straightforward computation yields
$$ (\JL \cdot \JM)=-\summ \j (L\cdot D_2)(M\cdot D_2),
$$

$$
\eqalign{ (\KL \cdot \KM)=-\summm &\k
\Bigl( (L\cdot E)\bigl( (M\cdot D_1) -(M\cdot D) \bigr) +\cr &-(L\cdot
D_2)\bigl(
3(M\cdot D_2)+(M\cdot E)\bigr) \Bigr),\cr }
$$

$$
\eqalign{ (\KKL \cdot \KKM)=-\summm &\kk
\Bigl( (L\cdot E)\bigl( 2(M\cdot D_1) -2(M\cdot D) \bigr) +\cr &-(L\cdot
D_2)\bigl(
3(M\cdot D_2)+2(M\cdot E)\bigr) \Bigr) \cr }
$$
and

$$
\eqalign{
 (\HL \cdot \HM)=  -{1\over 2}\summmm \h &\Bigl( -(L\cdot D)(M\cdot D) +
(L\cdot
D)(M\cdot D_1) \cr &+ (L\cdot D_1)(M\cdot D)-(L\cdot D_1)(M\cdot D_1) \cr
&-(L\cdot
D_2)(M\cdot D_2)- (L\cdot D_3)(M\cdot D_3) \Bigr) \cr}.
$$

In the last formula, we divide by $2$ because $D_2$ and $D_3$ are not
distinguished
from one another.

Now we are ready to write down the three formulas that we mentioned at the
beginning
of this chapter. Before we carry out the computation, we can explain briefly the
procedure. We have to look at the relation $(*)$ and  keep track of the Severi
degrees on which the characteristic numbers depend.

\

\ni (a) The first that relation we shall use is
$$ (\pull C \cdot \pull C) = 3 \ND
$$ This will give a formula expressing
$$
\ND \quad  {\rm in \ terms \ of} \quad \NDD \quad {\rm and}  \quad  \tDDD
\quad {\rm
with} \quad  D'< D \  {\rm and
\  } \  D''<D-E.
$$ This is clear; since $(C\cdot E)=0$ if we apply $(LM)$ to $L=M={\cal
O}_S(C)$ the
Severi degree $\fDE$ disappears.

\

\ni (b) Now we need a formula for $\tD$. We will imitate what we did to
compute the
degree of the variety of rational curves tangent to a fixed line in the
plane. We
define
$\E$ to be the class of the irreducible component of $\pi ^{-1}(E)$ that
dominates $B$.
Then we have
$$
\E = \pull E - 3{\cal G}_E
 -2 {\cal K}_E  - {\cal K}_0 -2{\cal K}'_E - {\cal K}'_0 -{\cal H}_E.
$$ This is obtained  as follows: for the coefficient  of ${\cal G}_E$ we  notice
that for any $b \in B_G$ we have
$(\E \cdot \Gb ) = 0$, while on the other hand $(\pull E \cdot \Gb ) = -3$.
The same
procedure yields the remaining terms.
$$
\tD =\E ^2 + (\E \cdot \dual ) .$$ This will give a recursion expressing
$$
\tD \  \ {\rm in \  terms\  of \  } \  \ND, \  \NDD ,\  \fDE ,\   {\rm and
\  }   \tDDD  \   \  {\rm
with
\  }\   D'< D \  {\rm and
\  } \  D''<D-E.
$$

\
\ni (c) The third and last step will be to find a formula for $\fDE$. This
will be
done by using
$$ (\pull F \cdot \pull F)=0
$$ which, as one can imply using $(*)$,  will give
$$
\fDE \  \  \  {\rm in \  terms\  of \  } \NDD \  \  {\rm and \  }  \  \tDDD
\   \
{\rm with \  }\   D'< D \  {\rm and
\  } \  D''<D-E.
$$

\

\ni
\underbar{Example.} If $D=2C$ we have on one hand that $(\pull C\cdot \pull
C)=3N(2C)$,   and on the other hand our formulas give
$$\eqalign{ (\pull C\cdot \pull C)=3[-12 &A^2 - 3 j(C,C) -k(C+2F,F) - 25
k'(F,C+2F) \cr
&\quad -14h(F,F,C+F)-2h(C+F,F,F)] \cr}
$$

\

Here are the relevant numbers  for the case $D=2C$.

$$ \eqalign{& j (C,C) = 105 \cr
&k(C+2F,F)=14 \cr
&k'(F,C+2F)=2 \cr
&h(F,F,C+F)=7 \cr
&  h(C+F,F,F)=21 \cr
&\nJ = 60 \cr
&\nK = \nKK = 2 \cr
&\nH = 13 \cr}$$
so that
$$A^2 =49$$
and we can conclude that $N(2C)=69$.

\

We will now state our main result for $\F3$:

\vfill\eject

\

\proclaim Theorem. Let $D\in \Pic (\F3 )$. Let $\ND $ be the number of
irreducible
rational curves in  $|D|$ that pass through $\rD$ general points.  Then

 $$
\eqalign{ &\ND =\cr &= {1\over 3}\summ N(D_1 )N(D_2)\Prod [{ \rD -3 \choose
\rDu -1}
(D_1 \cdot C)(D_2 \cdot C) - { \rD -3 \choose \rDu -2} (D_2 \cdot C)^2]+ \cr
 &+\summm N_2(D_1 )N(D_2) (E\cdot D_2) \Bigr[ { \rD -3 \choose \rDu -2}
(D_1 \cdot
C)(D_2 \cdot C) - { \rD -3 \choose \rDu -3} (D_2 \cdot C)^2\Bigl] \cr
 &+\summm N(D_1 )N_2(D_2) (E\cdot D_1) \Bigr[ { \rD -3 \choose \rDu -1}
(D_1 \cdot
C)(D_2 \cdot C) - { \rD -3 \choose \rDu -2} (D_2 \cdot C)^2\Bigl] \cr
 &+{1\over 3}\summmm  N(D_1 )N(D_2)N(D_3) (E\cdot D_1)(E\cdot D_2)(E\cdot
D_3)\cdot
\cr  & \cdot \Bigr[   { \rD -3 \choose \rDu -1 ,\  \  \rDt -1} [2(C\cdot
D_1)(C\cdot
D_2)+ (C\cdot D_1)(C\cdot D_3)+(C\cdot D_2)(C\cdot D_3)- (C\cdot D_3)^2 ] -\cr
 &  - { \rD -3 \choose \rDu -2 ,\  \  \rDt}[(C\cdot D_2)^2 + (C\cdot D_3)^2
+(C\cdot
D_2)(C\cdot D_3)]
\Bigl]
\cr}.
$$

\

\ni {\it Proof.} We just have to compute.  Applying (*) to $\pi^*C$ gives
$$
\eqalign{ (\pull C)^2 = &- (C\cdot D)^2A^2 - \summ \j (C\cdot D_2)^2  -\cr
&- \summm
3\k  (C\cdot D_2)^2 - \cr &- \summm 3\kk  (C\cdot D_2)^2 - \cr &- \summmm
\h [(C\cdot
D_2)^2+(C\cdot D_3)^2+ (C\cdot D_2)(C\cdot D_3)]\cr}
$$ This gives

 $$
\eqalign{ &\ND =\cr &= {1\over 3}\summ N(D_1 )N(D_2)\Prod [{ \rD -3 \choose
\rDu -1}
(D_1 \cdot C)(D_2 \cdot C) - { \rD -3 \choose \rDu -2} (D_2 \cdot C)^2]+ \cr
 &+\summm N_2(D_1 )N(D_2) (E\cdot D_2) [2{ \rD -3 \choose \rDu -2} (D_1
\cdot C)(D_2
\cdot C) -\cr
 & \  \  \  \  \  \  \  \  \  \  \  \  \  \  \  \  \  \   -{ \rD -3 \choose
\rDu -1}
(D_1 \cdot C)^2 -{ \rD -3 \choose \rDu -3} (D_2 \cdot C)^2]+    \cr
 &+{1\over 3}\summmm  N(D_1 )N(D_2)N(D_3) (E\cdot D_1)(E\cdot D_2)(E\cdot D_3)
\Bigr[  (C\cdot D)^2 { \rD -3 \choose \rDu -1 ,\  \  \rDt -1} - \cr
 & \  \  \  \  \  \  \  \  \  \  \  \  \  \  \  \  \  \
 -\bigl( (C\cdot D_2)^2 + (C\cdot D_3)^2 +(C\cdot D_2)(C\cdot D_3)\bigr)
 { \rD -2 \choose \rDu -1 ,\  \  \rDt } \Bigl]
\cr}.
$$ And this concludes the proof.

\

\

We will write as well the formulas for  the  degrees of the other loci that
we need. The
first formula is obtained by
 $$
\tD =\E ^2 + (\E \cdot \dual ),
$$ which gives

$$
\eqalign{
\tD = &-3\ND + 9 \fDE + (E\cdot D)A^2  +\cr  & + \summ \j (E \cdot D_2) +
\cr & +
\summm 6(\k + \kk ) + \cr & +\summmm \h [2(E\cdot D_2) +2(E\cdot D_3)-1] .\cr}
$$
Finally, the degree of the Severi variety parametrizing rational curves
having a point
of contact of order at least 3 with $E$ is obtained by $\pull F ^2=0$, which
translates into

$$
\eqalign{N_3(&D-E) =  \cr
&-(F\cdot D)^2A^2  - \cr & - \!\!\!\! \summ \j (F \cdot D_2)^2 - \cr & -
\!\!\!\! \summm \k [1+2(F\cdot D_2) + 3 (F\cdot D_2)^2 ]-\cr & - \!\!\!\!
\summm \kk
[2+4(F\cdot D_2) + 3 (F\cdot D_2)^2 ]- \cr & - \!\!\!\!\!\!\!\! \summmm \h
[(F\cdot D)^2
-2(F\cdot D)(F\cdot D_1) +(F\cdot D_1)^2 +(F\cdot D_2)^2+ (F\cdot D_3)^2]\cr}
$$

\

\

\cl{\bf References}

\

\ni [CH] L. Caporaso and  J. Harris, {\it Parameter spaces for curves on
surfaces and enumeration of rational
curves} \ 1995, Preprint.

\ni [CM]  B. Crauder and R. Miranda,
{\it  Quantum cohomology of rational surfaces.} The Moduli Space of Curves,
Progress
in mathematics 129,   Birkhauser (1995)  pp. 33-80.

\ni [DH] \ S. Diaz and J. Harris, {\it Geometry of the Severi variety}. \
Trans.A.M.S. 309 (1988),
pp. 1-34.

\ni [DI]  P. Di Francesco and C. Itzykson,
{\it  Quantum intersection rings.} The Moduli Space of Curves, Progress in
mathematics
129,   Birkhauser (1995)  pp. 81-148.

\ni [H] \  J. Harris, {\it On the Severi problem}.  \ Invent. Math. 84
(1986), pp. 445-461.

\ni [K] \ J. Koll\'ar, {\it Rational Curves on Algebraic Varieties},
Springer 1996.

\ni [KP] S. Kleiman and R. Piene,
Private communication.

\ni [KM] M. Kontsevich and Y. Manin,  \ {\it Gromov-Witten classes, quantum
cohomology
and enumerative geometry.}

\ni [P] R. Pandharipande, {\it Intersections of $\Q$-divisors on
Kontsevich's moduli space
$\overline{M}_{0,n}(\P^r, d)$ and enumerative geometry}. \ 1995,
Preprint.

\ni [R] Z. Ran, {\it Enumerative geometry of singular plane curves.}
Invent. Math. 97 (1989), pp. 447-465.

\ni [RT] Y. Ruan and G. Tian, {\it A
mathematical theory of quantum cohomology}, J. Diff. Geom. 42 No. 2 (1995)
pp. 259-367).

\end